\input amstex
\documentstyle{amsppt}
\magnification 1200
\vcorrection{-8mm}
\NoBlackBoxes
\input epsf

\def\Z{\Bbb Z}
\def\Q{\Bbb Q}

\def\CP{\Bbb CP}
\def\ord{\operatorname{ord}}
\def\Res{\operatorname{Res}}
\def\Sing{\operatorname{Sing}}
\def\Discr{\operatorname{Discr}}

\def\refAB        {1}
\def\refACC       {2}
\def\refDegtPacif {3}
\def\refDegtACC   {4}
\def\refDegtMax   {5}
\def\refHH        {6}
\def\refvH        {7}
\def\refOkaPho    {8}
\def\refOrevFA    {9}
\def\refOrevWebPg {10}
\def\refSW        {11}
\def\refSerre     {12}
\def\refYang      {13}

\def\eqHilbS      {1}

\rightheadtext{Parametric equations of plane sextic curves}

\topmatter
\title    Parametric equations of plane sextic curves
          with a maximal set of double points
\endtitle
\author   S.~Yu.~Orevkov
\endauthor
\abstract
We give explicit parametric equations for all irreducible
plane projective sextic curves which have at most double points
and whose total Milnor number is maximal (is equal to $19$).
In each case we find a parametrization over a number field of the
minimal possible degree and try to choose coordinates so that
the coefficients are as small as we can do.
\endabstract
\address IMT, Univ. Toulouse 3, France and Inst. Steklov, Moscow, Russia
\endaddress
\endtopmatter

\centerline{\it \hskip 40mm To the memory of Shriram Abhyankar } 

\bigskip

\document
Let $C$ be a sextic curve in $\CP^2$ with simple singularities. 
The total Milnor number $\mu(C)$ of $C$
(i.~e., the sum of Milnor numbers of singular points)
is at most $19$. Indeed, if $X$ is the double
covering of $\CP^2$ branched at $C$ and
$\sigma:\tilde X\to X$ is the resolution of singularities of $X$,
then $\tilde X$ is a K3 surface, 
hence the negative index of inertia of the
intersection form on $H_2(\tilde X)$ is $19$ and the
exceptional curves of $\sigma$ generate a
sublattice of rank $\mu(C)$ with a negative definite form.

Using a more detailed analysis of sublattices of $H_2(\tilde X)$, a complete
list of all maximal sets of simple singularities is obtained by
Yang [\refYang]. If $\mu(C)=19$, then $C$ is
rational by the genus formula and $C$ is rigid
because the K3 surface is rigid.
The latter means that there are only finitely many of
such curves up to projective equivalence.

As in [\refDegtACC], in this paper we consider {\bf irreducible 
sextic curves of total Milnor number 19 which have $A_n$
singularities only}.
We give explicit parametric equations for all of them.
There are 39 equisingularity classes of such curves. 
They are listed in Degtyarev's papers [\refDegtPacif--\refDegtMax]
where the number of projectively non-equivalent curves is computed
in each case.
For the equisingularity classes
(sets of singularities), we use
the same numbering as in [\refDegtPacif--\refDegtMax].

Explicit defining equations
are computed in [\refDegtACC, \refDegtMax]
for 33 of these 39 equisingularity classes (some of them were earlier
computed in [\refACC, \refDegtPacif, \refOkaPho, \refOrevFA]).
The computations in [\refDegtACC, \refDegtMax] are based on
the method proposed by Artal, Carmona, and Cogolludo in [\refACC]
(the ACC-method).

Here we use a more straightforward method which we applied already in
[\refOrevFA]. Namely, we just write a parametrization of $C$
with indeterminate coefficients and solve the simultaneous
equations imposed by the singularity types.
It happens that this method works better in the cases when the
most singular points have one branch and when it is clear how
to choose coordinates to cancel the action of $PSL(3)$ and $PSL(2)$.
Otherwise, for example, in the case of the $A_{19}$ singularity
treated in [\refACC], the ACC-method works better.

In the cases when the straightforward method does not work well, we
computed parametrizations by a ``reengineering''
from the defining equations from [\refDegtACC, \refDegtMax].

In all the 39 cases we find a parametrization over a field $E$
of the minimal possible degree. In all cases except four, we have
$E=F$ where $F$ is the minimal field of
definition of the curve. In the remaining four cases
(namely, {\bf 1}, {\bf 16}, {\bf 34}, {\bf 36}) $E$ is a quadratic extension of $F$.

By rephrasing Abhyankar, one could say:
{\it a nice curve must have a nice
equation}. For the curves in question I tried to find equations
as nice as I could (at least as short as I could).

I am grateful to Alexei Skorobogatov for useful discussions.

\head\S1. The parametric equations
\endhead
All the parametrizations below are available on the web page [\refOrevWebPg] in Maple format.
%
The condition $\mu(C)=19$ implies that exactly one singularity has two local
branches and all the others are irreducible. In the parametrizations below,
$p(t)$ is always a quadratic polynomial whose roots are mapped to
the reducible singularity (the $A_n$ with odd $n$). In the brackets we
list the values of $t$ mapped to the irreducible singularities in the order
they appear in the notation ``$A_{n_1}+A_{n_2}+\dots$''.

\smallskip
\noindent
{\bf 1.} $A_{19}$, see [\refACC].
 $F=\Q(\sqrt5)$, $E=F(\sqrt{-3})$, $\omega^2+\omega+1=0$,
$a = 3\pm\sqrt5$, $p=t$,
$$
\split
&x= \omega t^2 + a t^3 - \omega^2 t^4,\\
&y= 2(t+t^5) - 5a(t^2-t^4) + (6+3a)t^3,\\
&z=  (\omega^2-\omega t^6) + 4a(\omega t+\omega^2 t^5) +
 (81a/2-34)(t^2-t^4) + (63a-54)t^3.
\endsplit
$$
A parametrization over $F$ does not exist because the curve
does not have smooth real points.
\if01                                                              %
alias(a=RootOf(Z^2-6*Z+4)); alias(w=RootOf(Z^2+Z+1));              #
x := w*t^2 + a*t^3 - w^2*t^4;                                      #
y := 2*(t+t^5) - 5*a*(t^2-t^4) + (6+3*a)*t^3;                      #
z :=  (w^2-w*t^6) + 4*a*(w*t+w^2*t^5) +                            #
 (81*a/2-34)*(t^2-t^4) + (63*a-54)*t^3;                            #
F := factor(resultant(x*Z-z*X,y*Z-z*Y,t)/Z^6);                     #
F:=factor(F/op(1,F)/op(2,F));                                      #
factor(discrim(F,X));                                              #
factor(discrim(F,Y));                                              #
factor(discrim(F,Z));                                              #
F1:=subs(a=3+sqrt(5),F);                                           #
F2:=subs(a=3-sqrt(5),F);                                           #
d1:=factor(discrim(F1,Z)); d1:=factor(d1/op(1,d1)/op(2,d1));       #
fsolve(subs(X=1,d1),complex);                                      #
factor(subs(X=1,Y=0,F1)); fsolve(
d2:=factor(discrim(F2,Z)); d2:=factor(d2/op(1,d2)/op(2,d2));       #
fsolve(subs(X=1,d2),complex);                                      #
factor(subs(X=1,Y=0,F2)); fsolve(
\fi                                                                %

\smallskip\noindent
{\bf 2.} $A_{18}+A_1$, see [\refACC].
$F=E=\Q(a)$, $a^3-2a-2=0$, $[\infty]$,
$p=3t^2 - 12 t + 5 a^2 + 8 a + 17$,
$$
\split
  x= \;&(t+3-a)(t+3+a)\,p,\qquad
  y= \; t^2 + 3a^2+4a-1,\\
  z= \;&\big[9t^4 + 72t^3 - (52a^2+100a-232)t^2 - (204a^2+528a-420)t
\\
     &- 92a^2-326a+587\big]\,p.
\endsplit
$$
\if01                                                              %
alias(a=RootOf(Z^3-2*Z-2));                                        #
p := 3*t^2 - 12*t + 5*a^2 + 8*a + 17;                              #
x := (t + 3 - a)*(t + 3 + a)*p;                                    #
y := t^2 + 3*a^2+4*a-1;                                            #
z := (9*t^4 + 72*t^3 - (52*a^2+100*a-232)*t^2                      #
     - (204*a^2+528*a-420)*t - 92*a^2-326*a+587)*p;                #
F := factor(resultant(x*Z-z*X,y*Z-z*Y,t)/Z^6);                     #
F:=factor(F/op(1,F)/op(2,F));                                      #
factor(discrim(F,X));                                              #
factor(discrim(F,Y));                                              #
factor(discrim(F,Z));                                              #
\fi                                                                %

\smallskip\noindent
{\bf 3.} $(A_{17}+A_2)$, see [\refACC, \refOkaPho].
$E=F=\Q$, $[\infty]$, $p=t^2-3$,
$$
   x = (t^2 + 9)p^2,\qquad
   y = tp,\qquad
   z = t^4 - 12 t^2 + 3;
$$
the curve is symmetric under $(x:y:z)\mapsto(x:-y:z)$, $t\mapsto-t$.
\if01                                                              %
# x := (1+t^2)*(1-3*t^2)^2;                                        #
# y := t*(1-3*t^2);                                                #
# z := 9 - 76*t^2 + 51*t^4 + 72*t^6;                               #
x := (t^2 + 9)*(t^2 - 3)^2;                                        #
y := t*(t^2 - 3);                                                  #
z := t^4 - 12*t^2 + 3;                                             #
F := factor(resultant(x*Z-z*X,y*Z-z*Y,t)/Z^4);                     #
factor(discrim(F,X));                                              #
factor(discrim(F,Y));                                              #
factor(discrim(F,Z));                                              #
\fi                                                                %

\smallskip\noindent
{\bf 4.} $A_{16}+A_3$. see [\refACC, \refDegtMax].
$E=F=\Q(17)$, $a=(9\pm\sqrt{17})/8$, $[\infty]$, $p=t^2+t-a+2/3$,
$$
  x=(t^2-3\,t - 6+11 a)\,p,\quad
  y=t^2-a,\quad
  z=(t^2-6\,t+6+5 a)\,p^2.
$$
\if01                                                              %
a := (9+sqrt(17))/8;                                              #
p := t^2 + t - a+2/3;                                              #
x := (t^2 - 3*t - 6+11*a)*p;                                       #
y := t^2 - a;                                                      #
z := (t^2 - 6*t + 6+5*a)*p^2;                                      #
F := factor(resultant(x*Z-z*X,y*Z-z*Y,t)/Z^6);                     #
factor(discrim(F,X));                                              #
factor(discrim(F,Y));                                              #
factor(discrim(F,Z));                                              #
\fi                                                                %

\smallskip\noindent
{\bf 5.} $A_{16}+A_2+A_1$. see [\refACC, \refOrevFA]. $E=F=\Q(a)$,
$a^3-a^2+a-3=0$, $[\infty,0]$,

\noindent
$p = 2t^2 + (8+7a-7a^2)(2t+1)$,
$$
\split
  x &= \big[t^2 + (8+13a-13a^2)t -(1672-3049a+1261a^2)/14\big]\,t^2,\\
  y &= (2t-1+16a-13a^2)(2t+9+12a-15a^2)\, t^2 p,\\
  z &= \big(2t^2 + (a-a^2)(12t-9)\big)\, p.
\endsplit
$$
\if01                                                              %
alias(a=RootOf(Z^3-Z^2+Z-3));                                      #
p := 2*t^2 + (8+7*a-7*a^2)*(2*t+1);                                #
x := (t^2 + (8+13*a-13*a^2)*t -(1672-3049*a+1261*a^2)/14)*t^2;     #
y := (2*t-1+16*a-13*a^2)*(2*t+9+12*a-15*a^2)*t^2*p;                #
z := (2*t^2 + (a-a^2)*(12*t-9))*p;                                 #
F := factor(resultant(x*Z-z*X,y*Z-z*Y,t)/Z^4);                     #
dx:=factors(discrim(F,X))[2];                                      #
dy:=factors(discrim(F,Y))[2];                                      #
dz:=factors(discrim(F,Z))[2];                                      #
\fi                                                                %

\smallskip\noindent
{\bf 6.} $A_{15}+A_4$, see [\refACC, \refDegtMax].
$E=F=\Q(i)$. $[\infty]$, $p=15\,t^2 - 7 + 6i$,
%
%
$$
\split
  x &= (145\,t^2 - 580i\,t - 133-146i)p,\\
  y &= (1649\,t+724-465 i)(170\,t-7-164i), \\
  z &= (5\,t^2+10i\,t+1+2i)\,p^2.
\endsplit
$$
\if01                                                              %
p := 15*t^2 - 7 + 6*I;                                             #
x := (145*t^2 - 580*I*t - 133 - 146*I)*p;                          #
y := (1649*t + 724 - 465*I)*(170*t - 7 - 164*I);                   #
z := (5*t^2 + 10*I*t + 1 + 2*I)*p^2;                               #
F := factor(resultant(x*Z-z*X,y*Z-z*Y,t)/Z^6);                     #
factor(discrim(F,X));                                              #
factor(discrim(F,Y));                                              #
factor(discrim(F,Z));                                              #
\fi                                                                %

\smallskip\noindent
{\bf 7.} $A_{14}+A_4+A_1$,  see [\refACC, \refDegtMax].

\noindent
$E=F=\Q(a)$,
$2a^6 - 6a^5 + 10a^4 - 5a^3 - 10a^2 + 4a + 8=0$, $[\infty,0]$,

\noindent
$p=t^2 - 12\, t - (338 a^5 - 638 a^4 + 458 a^3 + 771 a^2 - 4382 a - 2592)/187$,
$$\split
x&=\big(t^2 + (4a-12)t - 2a^5 + 6a^4 - 2a^3 - 3a^2 - 14a + 24\big)\,t^2,\\
y&=\big(t^2 + 4at -(10a^5 - 22a^4 + 50a^3 - 73a^2 - 38a + 16)/11\big)\,p,\\
z&=\big(t^2 + (8a-8)t - (10a^5 - 22a^4 + 50a^3 - 73a^2 + 50a + 16)/11\big)\,p\,t^2.
\endsplit$$
\if01                                                              %
alias(a=RootOf(2*Z^6 - 6*Z^5 + 10*Z^4 - 5*Z^3 - 10*Z^2 + 4*Z + 8));#
p:=t^2 - 12*t - (338*a^5-638*a^4+458*a^3+771*a^2-4382*a-2592)/187; #
x:=(t^2 + (4*a-12)*t - 2*a^5+6*a^4-2*a^3-3*a^2-14*a+24)*t^2;       #
y:=(t^2 + 4*a*t -(10*a^5-22*a^4+50*a^3-73*a^2-38*a+16)/11)*p;      #
z:=(t^2+8*(a-1)*t-(10*a^5-22*a^4+50*a^3-73*a^2+50*a+16)/11)*p*t^2; #
F := factor(resultant(x*Z-z*X,y*Z-z*Y,t)/Z^6);                     #
factor(discrim(F,X));                                              #
factor(discrim(F,Y));                                              #
factor(discrim(F,Z));                                              #
\fi                                                                %
%
Note that $F$ is a cubic extension of $\Q(\sqrt{-15})$: we have $a=-{6\over11}b^3+b-{2\over11}$
where $b^3 = -{5\over4} + {11\over36}\sqrt{-15}$.

\smallskip\noindent
{\bf 8.} $(A_{14}+A_2)+A_3$, see [\refOkaPho]. $E=F=\Q$,
%
$[\infty,2]$,
 $p=t^2 + 2t - 9$,
$$
  x=t^2 - 5,\qquad
  y=(t^2 + 2t - 5)\,p, \qquad
  z=(t^2 + 4t - 11)\,p^2.
$$
\if01                                                              %
p := t^2 + 2*t - 9;                                                #
x := t^2 - 5;                                                      #
y := (t^2 + 2*t - 5)*p;                                            #
z := (t^2 + 4*t - 11)*p^2;                                         #
F := factor(resultant(x*Z-z*X,y*Z-z*Y,t)/Z^6);                     #
factor(discrim(F,X));                                              #
factor(discrim(F,Y));                                              #
factor(discrim(F,Z));                                              #
\fi                                                                %

\smallskip\noindent
{\bf 9.} $(A_{14}+A_2)+A_2+A_1$, see [\refOkaPho].
$E=F=\Q$, 
$[\infty,0,1]$,
$p = t^2 + 5\,t - 5$,
$$
x = (2t^2-3)\,p,\qquad
y = (t^2+5\,t+1)\,t^2,\qquad
z = (t^2 + 4t - 6)\,p\,t^2.
$$
\if01                                                              %
p := t^2 + 5*t - 5;                                                #
x := (2*t^2 - 3)*p;                                                #
y := (t^2 + 5*t + 1)*t^2;                                          #
z := (t^2+4*t-6)*p*t^2;                                            #
F := factor(resultant(x*Z-z*X,y*Z-z*Y,t)/Z^6);                     #
factor(discrim(F,X));                                              #
factor(discrim(F,Y));                                              #
factor(discrim(F,Z));                                              #
\fi                                                                %

\smallskip\noindent
{\bf 10.} $A_{13}+A_6$, see [\refDegtMax].
$E=F=\Q(a)$, $a^4+7=0$, $[\infty]$,
%
$p=3t^2 + 7a^3-5a^2-13a+39$,
$$\split
x&= (3t - 2a^3+2a^2-6a+2)(3t + 2a^3-a^2-2a+3)\,p,\\
y&= 3t^2 + 4(a^2-8a+5)\,t + (-1037a^3+2143a^2-1105a-2757)/33,\\
z&= \big[ 3t^2 - 2(a^2-8a+5)\,t + (31a^3+15a^2-149a+235)/3
\big]\,p^2.
\endsplit$$
\if01                                                              %
alias(a=RootOf(Z^4+7));                                            #
p := 3*t^2 + 7*a^3-5*a^2-13*a+39;                                  #
x := (3*t - 2*a^3+2*a^2-6*a+2)*(3*t + 2*a^3-a^2-2*a+3)*p;          #
y := 3*t^2+4*(a^2-8*a+5)*t+(-1037*a^3+2143*a^2-1105*a-2757)/33;    #
z := ( 3*t^2 - (a^2-8*a+5)*t + (31*a^3+15*a^2-149*a+235)/3 )*p^2;  #
F := factor(resultant(x*Z-z*X,y*Z-z*Y,t)/Z^4);                     #
factor(discrim(F,X));                                              #
factor(discrim(F,Y));                                              #
factor(discrim(F,Z));                                              #
\fi                                                                %

\smallskip\noindent
{\bf 11.} $A_{13}+A_4+A_2$, see [\refDegtMax].
$E=F=\Q(\sqrt{21})$,
$a=\pm\sqrt{21}$,
$[\infty,3-(3a/2)]$,
$p = 2t^2+11a-39$,
$q = 2t+3a-6$,
$$\split
x &= (t+a-3)(t-2a-3)\,p,\\
y &= (2t-5a-9)(2t-3a+9)\,q^2,\\
z &= (t-a+3)(t+3)\,p\,q^2.
\endsplit$$
\if01                                                              %
alias(a=RootOf(Z^2-21));                                           #
p := 2*t^2+11*a-39;                                                #
q := 2*t+3*a-6;                                                    #
x := (t+a-3)*(t-3-2*a)*p;                                          #
y := (2*t-5*a-9)*(2*t-3*a+9)*q^2;                                  #
z := (t-a+3)*(t+3)*p*q^2;                                          #
F := factor(resultant(x*Z-z*X,y*Z-z*Y,t)/Z^6);                     #
factor(discrim(F,X));                                              #
factor(discrim(F,Y));                                              #
factor(discrim(F,Z));                                              #
\fi                                                                %

\smallskip\noindent
{\bf 12.} $A_{12}+A_7$, see [\refDegtMax].
$E=F=\Q(i)$,
%
%
$[\infty]$,
$p=3t^2+1$,
$$\split
x &= \big(13\,t^2 + (-8+12i)\,t + 7-4i\big)\,p,\\
y &= \big(13\,t^2 + (-20+30i)\,t -11-16i\big)\,p^2,\\
z &= (13\,t + 3+2i)(13t + 1-8i).
\endsplit$$
\if01                                                              %
p := 3*t^2 + 1;                                                    #
x := (13*t^2 + (-8+12*I)*t + 7-4*I)*p;                             #
y := (13*t^2 + (-20+30*I)*t -11-16*I)*p^2;                         #
z := (13*t + 3+2*I)*(13*t + 1-8*I);                                #
F := factor(resultant(x*Z-z*X,y*Z-z*Y,t)/Z^4);                     #
factor(discrim(F,X));                                              #
factor(discrim(F,Y));                                              #
factor(discrim(F,Z));                                              #
\fi                                                                %

\smallskip\noindent
{\bf 13.} $A_{12}+A_6+A_1$, see [\refDegtMax].
$E=F=\Q(a)$,
$7a^3-16a^2+12a-4=0$,
$[\infty,0]$,
$p=3t^2 + (12a-18)t + 28a^2 - 53a + 36$,
$$\split
x &= \big(t^2+(4a-2)t+14a^2-15a+4\big)\,t^2,\\
y &= (t^2 + 4t - 7a + 4)\,p,\\
z &= \big(t^2 + (4a+4)t + 13a\big)\,p\,t^2.
\endsplit$$

\if01                                                              %
alias(a=RootOf(7*Z^3-16*Z^2+12*Z-4));                              #
p := 3*t^2 + (12*a-18)*t + 28*a^2 - 53*a + 36;                     #
x := (t^2+(4*a-2)*t+14*a^2-15*a+4)*t^2;                            #
y := (t^2 + 4*t - 7*a + 4)*p;                                      #
z := (t^2 + (4*a+4)*t + 13*a)*p*t^2;                               #
F := factor(resultant(x*Z-z*X,y*Z-z*Y,t)/Z^6);                     #
factor(discrim(F,X));                                              #
factor(discrim(F,Y));                                              #
factor(discrim(F,Z));                                              #
\fi                                                                %

\smallskip\noindent
{\bf 14.} $A_{12}+A_4+A_3$, see [\refDegtMax]. $E=F=\Q$,
$[\infty,0]$,
$p=85 t^2 - 442 t + 507$,
$$
x = (5t^2-12t-13)\,t^2,\qquad
y = (25t^2-20t-221)\,t^4,\qquad
z = t^2-4t+3.
$$
\if01                                                              %
x := (5*t^2-12*t-13)*t^2;                                          #
y := (25*t^2-20*t-221)*t^4;                                        #
z := t^2-4*t+3;                                                    #
F := factor(resultant(x*Z-z*X,y*Z-z*Y,t)/Z^4);                     #
factor(discrim(F,X));                                              #
factor(discrim(F,Y));                                              #
factor(discrim(F,Z));                                              #
\fi                                                                %

\smallskip\noindent
{\bf 15.} $A_{12}+A_4+A_2+A_1$.
$E=F=\Q(a)$, $15a^3-48a^2+40a-10=0$,

\noindent
$[\infty,0,30a^2 - {159\over 2}a + {55\over 2}]$,
$p=t^2+ (60a^2-141a+64)t + 3828/5\,a^2 - 1851\,a + 772$,
%
$$\split
x &= \big(t^2 + 9t -24a^2+24a+5\big)t^2,\\
y &= t^2 + 9at + 84a^2-165a+50,\\
z &= \big(t^2 + (-9a+18)t -51a^2+24a+41 \big)t^4.
\endsplit $$
\if01                                                              %
alias(a=RootOf(15*Z^3-48*Z^2+40*Z-10));                            #
x := (t^2 + 9*t -24*a^2+24*a+5)*t^2;                               #
y := t^2 + 9*a*t + 84*a^2-165*a+50;                                #
z := (t^2 + (-9*a+18)*t -51*a^2+24*a+41)*t^4;                      #
F := factor(resultant(x*Z-z*X,y*Z-z*Y,t)/Z^6);                     #
factor(discrim(F,X));                                              #
factor(discrim(F,Y));                                              #
factor(discrim(F,Z));                                              #
                                                                   #
restart;                                                           #
alias(b=RootOf(Z^3 - Z^2 - Z - 5));                                #
factor(15*Z^3-48*Z^2+40*Z-10,b);                                   #
a:=(2*b^2+3*b+13)/15;                                              #
simplify(subs(Z=a,15*Z^3-48*Z^2+40*Z-10));                         #
x := factor((t^2 + 9*t -24*a^2+24*a+5)*t^2);                       #
y := factor(t^2 + 9*a*t + 84*a^2-165*a+50);                        #
z := factor((t^2 + (-9*a+18)*t -51*a^2+24*a+41)*t^4);              #
\fi                                                                %

\smallskip\noindent
{\bf 16.} $A_{11}+2A_4$, see [\refDegtMax].
$F=\Q(\sqrt2)$, $E=F(i)$,
%
$a=\pm\sqrt2$, $c = 4a+34i-27$, 
$[(3+5i\pm\sqrt{30i\,}\,)/4]$,
$p=t^2+(2ai+a-3i-2)t+a+i-1$,
$$\split
x &= \big[17t^4 +(34ai-9a-34i-20)t^3 + (136ai+49a+34i-148)t^2\\
&\qquad\quad\;\; +(-42ai+127a+122i+48)t -37ai-13a-18i+24\big]\,p,\\
y &= c(t+2ai+a-3i-1)(17t-6ai-7a-15i-9)\big(2t^2-(5i+3)t-2\big)^2,\\
z &= (2t-i)\,p^2.
\endsplit$$
\if01                                                              %
a:=sqrt(2);  i:=I;                                                 #
c := 8*a*i-8*a-14*i+13;                                            #
c := 4*a+34*i-27;                                                  #
p := t^2+(2*a*i+a-3*i-2)*t+a+i-1;                                  #
x := (17*t^4 +(34*a*i-9*a-34*i-20)*t^3+(136*a*i+49*a+34*i-148)*t^2 #
  +(-42*a*i+127*a+122*i+48)*t -37*a*i-13*a-18*i+24)*p;             #
y := c*(t+2*a*i+a-3*i-1)*(17*t-6*a*i-7*a-15*i-9)*(2*t^2-(5*i+3)*t-2)^2;
z := (2*t-i)*p^2;                                                  #
F := factor(resultant(x*Z-z*X,y*Z-z*Y,t)/Z^6);                     #
factor(discrim(F,X));                                              #
factor(discrim(F,Y));                                              #
factor(discrim(F,Z));                                              #
\fi                                                                %
%
Another parametrization (with $0,\infty\mapsto A_{11}$) is: $E=\Q(b)$,
%
$b^2 = -7\pm4\sqrt2$, thus $b^4 + 14b^2 + 17 = 0$,
$$\split
x &= \big[ 32t^4-(7b^3-49b^2+17b-135)t^3 - (28b^3-44b^2+52b-4)t^2\\
&\qquad\quad\;\;+(10b^3+2b^2-26b-50)t +2b^3+4b^2+18b+32 \big]\,t,\\
y &= (b^2-2b+9)\big[32t^2+(28b^3+8b^2+28b+16)t+2b^3-b^2+24b+19\big]\\
&\qquad\quad\;\;\times\big[16t^2 - (b^3-7b^2+7b-17)t + b^3+b^2+15b-1\big]^2,\\
z &= (b^2+2b+9)\big(2t^2-(b^2-1)t-2\big)\,t^2;
\endsplit$$
\if01                                                               %
alias(b=RootOf(b^4 + 14*b^2 + 17));                                 #
x := (32*t^4-(7*b^3-49*b^2+17*b-135)*t^3-(28*b^3-44*b^2+52*b-4)*t^2 #
     +(10*b^3+2*b^2-26*b-50)*t + 2*b^3+4*b^2+18*b+32)*t;            #
y :=(b^2-2*b+9)*(32*t^2+(28*b^3+8*b^2+28*b+16)*t+2*b^3-b^2+24*b+19) #
   *(16*t^2 - (b^3-7*b^2+7*b-17)*t + b^3+b^2+15*b-1)^2;             #
z := (b^2+2*b+9)*(2*t^2-(b^2-1)*t-2)*t^2;                           #
F := factor(resultant(x*Z-z*X,y*Z-z*Y,t)/Z^5);                      #
factor(discrim(F,X));                                               #
factor(discrim(F,Y));                                               #
factor(discrim(F,Z));                                               #
\fi                                                                 %
A parametrization over $F$ does not exist because the curve
does not have smooth real points.

\smallskip\noindent
{\bf 17.} $(A_{11}+2A_2)+A_4$, see [\refOkaPho]. $E=F=\Q$,
%
$[\pm i\sqrt{2/3},\infty]$,
$p=3t^2+9t+8$,
$$
x = (3t^2-10)\,p^2,\qquad
y = (3t^2+9t+5)\,p,\qquad
z = 3t^2+3t+2.
$$
\if01                                                              %
p := 3*t^2+9*t+8;                                                  #
x := (3*t^2-10)*p^2;                                               #
y := (3*t^2+9*t+5)*p;                                              #
z := 3*t^2+3*t+2;                                                  #
F := factor(resultant(x*Z-z*X,y*Z-z*Y,t)/Z^4);                     #
factor(discrim(F,X));                                              #
factor(discrim(F,Y));                                              #
factor(discrim(F,Z));                                              #
dX:=factor(diff(x/z,t)*z^2);
dY:=factor(diff(y/z,t)*z^2);
\fi                                                                %

\smallskip\noindent
{\bf 18.} $A_{10}+A_9$, see [\refDegtMax]. $E=F=\Q(\sqrt5)$,
$a=\pm\sqrt5$, $[\infty]$, $p = t^2 - 11 - 22/3\,a$,
$$\split
x &= (t-3-2a)(t-5+4a)\,p,\\
y &= \big(t^2 - (32-8a)t + 149+70a\big)\,p^2,\\
z &= t^2 + (16-4a)t+161-226/5\,a.
\endsplit$$
\if01                                                              %
alias(a=RootOf(Z^2-5));                                            #
p := t^2 - 11 - 22/3*a;                                            #
x := (t - 3 - 2*a)*(t - 5 + 4*a)*p;                                #
y := (t^2 - (32-8*a)*t + 149 + 70*a)*p^2;                          #
z :=  t^2 + (16-4*a)*t + 161 - 226/5*a;                            #
F := factor(resultant(x*Z-z*X,y*Z-z*Y,t)/Z^4);                     #
factor(discrim(F,X));                                              #
factor(discrim(F,Y));                                              #
factor(discrim(F,Z));                                              #
\fi                                                                %

\smallskip\noindent
{\bf 19.} $A_{10}+A_8+A_1$, see [\refDegtMax]. $E=F=\Q(a)$,
$a^3-4a^2+8a-4=0$, $[\infty,0]$,

\noindent
$p=t^2 + (8a-6)t + 9a^2-29a+16$,
$$\split
x &= \big( t^2 + (4a-2)t + 15a^2-45a+24 \big)\,t^2,\\
y &= \big( t^2 - (4a-4)t + 5a^2-9a+4 \big)\,p,\\
z &= \big( t^2 + 4t - 3a^2 + 3a \big)\,p\,t^2.
\endsplit$$
\if01                                                              %
alias(a=RootOf(Z^3-4*Z^2+8*Z-4));                                  #
p := t^2 + (8*a-6)*t + 9*a^2-29*a+16;                              #
x := ( t^2 + (4*a-2)*t + 15*a^2-45*a+24 )*t^2;                     #
y := ( t^2 - (4*a-4)*t + 5*a^2-9*a+4 )*p;                          #
z := ( t^2 + 4*t - 3*a^2 + 3*a )*p*t^2;                            #
F := factor(resultant(x*Z-z*X,y*Z-z*Y,t)/Z^6);                     #
factor(discrim(F,X));                                              #
factor(discrim(F,Y));                                              #
factor(discrim(F,Z));                                              #
\fi                                                                %

\smallskip\noindent
{\bf 20.} $A_{10}+A_7+A_2$, see [\refDegtMax]. $E=F=\Q(\sqrt3)$,
%
$a=\pm\sqrt3$,
$[\infty,(3+27a)/22]$,
$p=2t^2-6+a$,
$$\split
x &= \big(22 t^2 - (12 + 20 a)t + 30 + 39 a \big)\,p,\\
y &= \big(22 t^2 - (42 + 26 a)t - 6 + 45 a \big)\,p^2,\\
z &= 22 t^2 + (18 - 14 a)t - (300 + 159 a)/11.
\endsplit$$
\if01                                                              %
alias(a=RootOf(Z^2-3));  p := 2*t^2 - 6 + a;                       #
x := (22*t^2 - (12 + 20*a)*t + 30 + 39*a)*p;                       #
y := (22*t^2 - (42 + 26*a)*t - 6 + 45*a)*p^2;                      #
z :=  22*t^2 + (18 - 14*a)*t - (300 + 159*a)/11;                   #
F := factor(resultant(x*Z-z*X,y*Z-z*Y,t)/Z^4);                     #
factor(discrim(F,X));                                              #
factor(discrim(F,Y));                                              #
factor(discrim(F,Z));                                              #
\fi                                                                %

\smallskip\noindent
{\bf 21.} $A_{10}+A_6+A_3$, see [\refDegtMax].
$E=F=\Q(\sqrt{-7})$,
%
%
$a=\pm i\sqrt7$,
$[0,\infty]$,

\noindent
 $p=22\,t^2+(27-3a)t+9-3a$,
$$\split
x &= (4t+1-a)(t+1)t^2,\\
y &= 22\,t^2+(56-16a)t-3-7a,\\
z &= \big(22\,t^2+(-1+5a)t-5+3a\big)\,t^4.
\endsplit$$
\if01                                                              %
alias(a=RootOf(Z^2+7));                                            #
x := (4*t+1-a)*(t+1)*t^2;                                          #
y := 22*t^2+(56-16*a)*t-3-7*a;                                     #
z := (22*t^2+(-1+5*a)*t-5+3*a)*t^4;                                #
F := factor(resultant(x*Z-z*X,y*Z-z*Y,t)/Z^6);                     #
factor(discrim(F,X));                                              #
factor(discrim(F,Y));                                              #
factor(discrim(F,Z));                                              #
\fi                                                                %

\smallskip\noindent
{\bf 22.} $A_{10}+A_6+A_2+A_1$. $E=F=\Q(a)$,
%
%
$a^3-a^2+3=0$,
$[\infty,0,4a^2-3]$,
$p=t^2 + (2a^2+6a+3)t - 26a^2+78a+123$,
$$\split
x &= \big( t^2 + 9t -22a^2+30a+21 \big)\,t^2,\\
y &= t^2 - (6a^2-6a+9)t + 50a^2 - 78a + 75,\\
z &= \big( t^2 + (6a^2-6a+27)t -22a^2+66a-33 \big)\,t^4.
\endsplit$$
\if01                                                              %
alias(a=RootOf(Z^3-Z^2+3));                                        #
x := ( t^2 + 9*t - 22*a^2+30*a+21 )*t^2;                           #
y :=   t^2 - (6*a^2-6*a+9)*t + 50*a^2 - 78*a + 75;                 #
z := ( t^2 + (6*a^2-6*a+27)*t -22*a^2+66*a-33 )*t^4;               #
F := factor(resultant(x*Z-z*X,y*Z-z*Y,t)/Z^6);                     #
factor(discrim(F,X));                                              #
factor(discrim(F,Y));                                              #
factor(discrim(F,Z));                                              #
\fi                                                                %

\smallskip\noindent
{\bf 23.} $A_{10}+A_5+A_4$, see [\refDegtMax]. $E=F=\Q(\sqrt{15})$,
$a=\pm\sqrt{15}$,
$[\infty,0]$,

\noindent
$p=5t^2 + (3a+5)t + (70-24a)/17$,
$$\split
x &= (t+a-4)(t+1)\,p,\\
y &= \big(t^2 + (4a-12)t - 2a+8\big)\,p\,t^2,\\
z &= \big(5t^2 + (8a-10)t + (110-2a)/7 \big)\,t^2.
\endsplit$$
\if01                                                              %
alias(a=RootOf(Z^2-15));                                           #
p := 5*t^2 + (3*a+5)*t + (70-24*a)/17;                             #
x := (t+a-4)*(t+1)*p;                                              #
y := (t^2 + (4*a-12)*t - 2*a+8)*p*t^2;                             #
z := (5*t^2 + (8*a-10)*t + (110-2*a)/7 )*t^2;                      #
F := factor(resultant(x*Z-z*X,y*Z-z*Y,t)/Z^4);                     #
factor(discrim(F,X));                                              #
factor(discrim(F,Y));                                              #
factor(discrim(F,Z));                                              #
\fi                                                                %

\smallskip\noindent
{\bf 24.} $A_{10}+2A_4+A_1$, see [\refDegtMax]. $E=F=\Q(a)$,
%
%
$a^3-a^2-a-1=0$,
$[\infty,\pm\sqrt{5a^2-20}]$,
$p=t^2 + (-12a^2+18a+8)t - 39a^2+12a+116$,
$$\split
x &= (t - a)(t + 4a^2-5a-4)\,p,\\
y &= \big[t^2 + (-20a^2+30a+8)t + 45a^2+20a-140\big](t^2 - 5a^2+20)^2,\\
z &= t^2 + (4a^2-6a)t + a^2+4a+4.
\endsplit$$
\if01                                                              %
alias(a=RootOf(Z^3-Z^2-Z-1));                                      #
x := (t - a)*(t + 4*a^2-5*a-4)*                                    #
          ( t^2 + (-12*a^2+18*a+8)*t - 39*a^2+12*a+116 );          #
y := (t^2 + (-20*a^2+30*a+8)*t +45*a^2+20*a-140)*(t^2-5*a^2+20)^2; #
z := t^2 + (4*a^2-6*a)*t + a^2+4*a+4;                              #
F := factor(resultant(x*Z-z*X,y*Z-z*Y,t)/Z^4);                     #
factor(discrim(F,X));                                              #
factor(discrim(F,Y));                                              #
factor(discrim(F,Z));                                              #
\fi                                                                %

\smallskip\noindent
{\bf 25.} $A_{10}+A_4+A_3+A_2$, see [\refDegtMax]. $E=F=\Q$,
%
%
$[\infty,0,4]$,
$p=20t^2-55t-121$,
$$
x = (7t^2 - 35t + 22)t^2,\qquad
y = (8t^2 - 52t + 77)t^4,\qquad
z = 2t^2 - 7t - 7.
$$
%
\if01                                                              %
x := (7*t^2 - 35*t + 22)*t^2;                                      #
y := (8*t^2 - 52*t + 77)*t^4;                                      #
z := 2*t^2 - 7*t - 7;                                              #
F := factor(resultant(x*Z-z*X,y*Z-z*Y,t)/Z^4);                     #
factor(discrim(F,X));                                              #
factor(discrim(F,Y));                                              #
factor(discrim(F,Z));                                              #
\fi                                                                %

\smallskip\noindent
{\bf 26.} $A_{10}+A_4+2A_2+A_1$. $E=F=\Q(\sqrt5)$,
%
%
$a=\pm\sqrt5$, $b=a-1$,

\noindent
$[\infty,0$,
roots of $t^2 + 3t +(15-a)/6]$,
$p=t^2+(a-4)t+(15-23a)/6$,
$$
x = \big(6t^2 + 6t + 9a - 25\big)t^2,\quad
y =   6t^2 + 6at + 5b,               \quad
z = \big( 6t^2 + (12-6a)t - 11b \big)\,t^4.
$$
\if01                                                              %
alias(a=RootOf(Z^2-5)); b:=a-1;                                    #
x := (6*t^2 + 6*t + 9*a-25)*t^2;                                   #
y :=  6*t^2 + 6*a*t + 5*b;                                         #
z := (6*t^2 + (12-6*a)*t - 11*b)*t^4;                              #
F := factor(resultant(x*Z-z*X,y*Z-z*Y,t)/Z^6);                     #
factor(discrim(F,X));                                              #
factor(discrim(F,Y));                                              #
factor(discrim(F,Z));                                              #
\fi                                                                %

\smallskip\noindent
{\bf 27.} $A_9+A_6+A_4$, see [\refDegtMax]. $E=F=\Q(a)$,
$a^3-5a-5=0$,
$[0,\infty]$,

\noindent
$p = 3t^2 - (4a^2+4a-7)t + 4a^2 + 13a + 20$,
$$\split
x &= \big[ 3t^2 - (4a^2+4a-25)t - (274a^2-365a-1510)/31 \big]\,t^2,\\
y &= (t - a^2 - a + 1)(3t - a^2 - a - 5)\,t^2 p,\\
z &= \big[ t^2 + 6t - 2a^2 - 11a - 10 \big]\,p.
\endsplit$$
\if01                                                              %
alias(a=RootOf(Z^3-5*Z-5));                                        #
p := 3*t^2 - (4*a^2+4*a-7)*t + 4*a^2 + 13*a + 20;                  #
x := ( 3*t^2 - (4*a^2+4*a-25)*t - (274*a^2-365*a-1510)/31 )*t^2;   #
y := (t - a^2 - a + 1)*(3*t - a^2 - a - 5)*t^2*p;                  #
z := (t^2 + 6*t - 2*a^2 - 11*a - 10)*p;                            #
F := factor(resultant(x*Z-z*X,y*Z-z*Y,t)/Z^4);                     #
factor(discrim(F,X));                                              #
factor(discrim(F,Y));                                              #
factor(discrim(F,Z));                                              #
\fi                                                                %

\smallskip\noindent
{\bf 28.} $A_9+2A_4+A_2$, see [\refDegtMax]. $E=F=\Q$,
%
$[\pm i\sqrt{15/2},\infty]$,
$p=2t^2-5$,
$$
x = pt,\qquad
y = 4t^4 - 80t^2 + 15,\qquad
z = (2t^2 + 75)p^2;
$$
the curve is symmetric under $(x:y:z)\mapsto (-x:y:z)$, $t\mapsto-t$.
\if01                                                              %
p := 2*t^2 - 5;                                                    #
x := p*t;                                                          #
y := 4*t^4 - 80*t^2 + 15;                                          #
z := (2*t^2 + 75)*p^2;                                             #
F := factor(resultant(x*Z-z*X,y*Z-z*Y,t)/Z^6);                     #
factor(discrim(F,X));                                              #
factor(discrim(F,Y));                                              #
factor(discrim(F,Z));                                              #
\fi                                                                %

\smallskip\noindent
{\bf 29.} $(2A_8)+A_3$, see [\refOkaPho]. $E=F=\Q$, 
$[0,\infty]$,
$p=t^2+2t-1$,
$$
x = 3t^2-4t+1,\qquad
y = (3t^2+4t-3)t^2,\qquad
z = (t^2+4t+3)t^4;
$$
the curve is symmetric under $(x:y:z)\mapsto(z:-y:x)$, $t\mapsto-1/t$.
\if01                                                              %
x := 3*t^2-4*t+1;                                                  #
y := (3*t^2+4*t-3)*t^2;                                            #
z := (t^2+4*t+3)*t^4;                                              #
F := factor(resultant(x*Z-z*X,y*Z-z*Y,t)/Z^6);                     #
factor(discrim(F,X));                                              #
factor(discrim(F,Y));                                              #
factor(discrim(F,Z));                                              #
\fi                                                                %

\smallskip\noindent
{\bf 30.} $A_8+A_7+A_4$, see [\refDegtMax]. $E=F=\Q(i)$,
$[\infty,(i-1)/2]$,
%
$p = 34t^2 - 8+19i$,
$$\split
x &= \big(30\,t^2+20(1-i)t-4-3i\big)\,p,\\
y &= \big(10t^2-(2+6i)t+5i\big)\,p^2,\\
z &= (10t+9-3i)(30t+19-13i).
\endsplit$$
\if01                                                              %
p := 34*t^2 - 8+19*I;                                              #
x := (30*t^2+20*(1-I)*t-4-3*I)*p;                                  #
y := (10*t^2-(2+6*I)*t+5*I)*p^2;                                   #
z := (10*t+9-3*I)*(30*t+19-13*I);                                  #
F := factor(resultant(x*Z-z*X,y*Z-z*Y,t)/Z^4);                     #
factor(discrim(F,X));                                              #
factor(discrim(F,Y));                                              #
factor(discrim(F,Z));                                              #
\fi                                                                %

\smallskip\noindent
{\bf 31.} $A_8+A_6+A_4+A_1$, see [\refDegtMax]. $E=F=\Q(a)$,
%
%
%
$a^3 - a^2 - a + 5 = 0$,
$[\infty,0,1]$,

\noindent
$p=t^2+(a^2-4a+5)t + (43a^2-104a+120)/7$,
$$\split
  x &= \big[ 15t^2 + (5-20a+5a^2)t - 70 + 64a - 23a^2]\,t^2,\\
  y &= \big[ 21t^2 + (49 - 28a + 7a^2)t -18 + 24a - 11a^2 \big]\,(t-1)^2,\\
  z &= \big[ 5t^2 + (1-4a-a^2)t - 10+5a^2 \big]\,t^2\,(t-1)^2.
\endsplit$$
\if01                                                              %
alias(a=RootOf(Z^3 - Z^2 - Z + 5));                                #
x := (15*t^2 + (5-20*a+5*a^2)*t - 70 + 64*a - 23*a^2)*t^2;         #
y := (21*t^2 + (49-28*a+7*a^2)*t -18 + 24*a - 11*a^2)*(t-1)^2;     #
z := (5*t^2 + (1-4*a-a^2)*t - 10+5*a^2)*t^2*(t-1)^2;               #
F := factor(resultant(x*Z-z*X,y*Z-z*Y,t)/Z^6);                     #
factor(discrim(F,X));                                              #
factor(discrim(F,Y));                                              #
factor(discrim(F,Z));                                              #
\fi                                                                %

\smallskip\noindent
{\bf 32.} $(A_8+A_5+A_2)+A_4$, see [\refOkaPho]. $E=F=\Q(i)$,
$[\infty,1,0]$,
$p=t^2+2it-(8+4i)/5$,
$$
    x = (t^2-4)pt^2,\quad
    y = (3t^2-2t+2)p, \quad
    z=(3t^2 - (2-6i)t - 4-4i)t^2.
$$
\if01                                                              %
p := t^2+2*I*t-(8+4*I)/5;                                          #
x := (t^2 - 4)*p*t^2;                                              #
y := (3*t^2 - 2*t + 2 )*p;                                         #
z := (3*t^2 - (2-6*I)*t - 4-4*I)*t^2;                              #
F := factor(resultant(x*Z-z*X,y*Z-z*Y,t)/Z^4);                     #
factor(discrim(F,X));                                              #
factor(discrim(F,Y));                                              #
factor(discrim(F,Z));                                              #
\fi                                                                %

\smallskip\noindent
{\bf 33.} $(A_8+3A_2)+A_4+A_1$, see [\refOkaPho]. $E=F=\Q$,
$[\infty$, roots of $t^3-3t-3,\, 0]$,
$p=t^2+3t-9$,
$$
x = (t^2+t-3)t^2,\qquad
y = (t^2+3t+3)t^4,\qquad
z = t^2-t-1.
$$
\if01                                                              %
x := (t^2+t-3)*t^2;                                                #
y := (t^2 + 3*t + 3)*t^4;                                          #
z := t^2 - t - 1;                                                  #
F := factor(resultant(x*Z-z*X,y*Z-z*Y,t)/Z^4);                     #
factor(discrim(F,X));                                              #
factor(discrim(F,Y));                                              #
factor(discrim(F,Z));                                              #
\fi                                                                %

\smallskip\noindent
{\bf 34.} $A_7+2A_6$, see [\refDegtMax]. $F=\Q(\sqrt{-7})$, $E=F(\sqrt{-3})$,
%
%
$a=\pm i\sqrt7$, $b=\pm i\sqrt3$, $p=t$,
$$\split
x & =4(a+1)\big[16t^4-(ab+a+b-23)t^3-(12ab-20a-28b-60)t^2\\
  & \qquad\quad\;\; -(23ab+2a+3b-106)t-10ab+46\big]t,\\
y & = (ab+5)\big[8t^2 +(2ab-16a-26b-12)t +5ab+15a+23b-23\big]\\
  & \qquad\quad\;\; \times\big[8t^2 -(ab+3a+5b-9)t-ab-3a-5b+5\big]^2,\\
z & = 2(a+1)^2(ab-5)\big(2t^2-(a-3)t+1\big)\,t^2.
\endsplit$$
The implicit equation has coefficients in $F$:
$$
\split
&1/4\,y^4\;+\;\big[x^2+(2-6a)x+58+(23/2)a\big]y^3\\
&\; +\big[x^4 -(10+26a)x^3 + (162-99a)x^2 + (1883-265a)x + 1799+1428a\big]y^2\\
&\; -(a+1)\big[7x^3 + (79-2a)x^2 + (182+9a)x + 658 + 336a\big]h(x)y\\
&\; -\big[(15+7a)x^2 + (175+43a)x +623-21a\big]h(x)^2,
\qquad h=4x^2+4ax-7+7a,
\endsplit
$$
and the $A_6$ points are on the line $y=0$ at the roots of $h$.
\if01                                                             %
alias(a=RootOf(Z^2+7)); alias(b=RootOf(Z^2+3));                   #
x:=4*(a+1)*(16*t^4-(a*b+a+b-23)*t^3-(12*a*b-20*a-28*b-60)*t^2     #
          -(23*a*b+2*a+3*b-106)*t-10*a*b+46)*t;                   #
y:=(a*b+5)*(8*t^2 +(2*a*b-16*a-26*b-12)*t +5*a*b+15*a+23*b-23)    #
          *(8*t^2 -(a*b+3*a+5*b-9)*t-a*b-3*a-5*b+5)^2;            #
z:=2*(a+1)^2*(a*b-5)*(2*t^2-(a-3)*t+1)*t^2;                       #
F := factor(resultant(x*Z-z*X,y*Z-z*Y,t)/Z^5):                    #
F := F/op(1,F)/op(2,F);                                           #
factor(discrim(F,X));                                             #
factor(discrim(F,Y));                                             #
factor(discrim(F,Z));                                             #
h := 4*X^2+4*a*X-7+7*a;                                           #
f := 1/4*Y^4 + (X^2+(2-6*a)*X+58+(23/2)*a)*Y^3
 +(X^4 -(10+26*a)*X^3 + (162-99*a)*X^2 + (1883-265*a)*X + 1799+1428*a)*Y^2
 -(a+1)*(7*X^3 + (79-2*a)*X^2 + (182+9*a)*X + 658 + 336*a)*h*Y
 -((15+7*a)*X^2 + (175+43*a)*X +623-21*a)*h^2;                    #
factor(subs(Z=1,F/f));                                            #
\fi                                                               %

\smallskip\noindent
{\bf 35.} $A_7+A_6+A_4+A_2$, see [\refDegtMax].
$E=F=\Q(\sqrt{21})$,
$a=\pm\sqrt{21}$,
$[\infty, 0, 4]$,

\noindent
$p=2t^2-(9-a)t+(2a-78)/5$,
$$\split
x &= \big(2t^2-(23-3a)t+116-24a\big)\,t^2 p,\\
y &= (t^2-5t-2+4/3\,a)t^2,\\
z &= (2t^2-(1+a)t+4-8/3\,a)\,p.
\endsplit$$
\if01                                                              %
alias(a=RootOf(Z^2-21));                                           #
p := 2*t^2-(9-a)*t+(2*a-78)/5;                                     #
x := (2*t^2-(23-3*a)*t+116-24*a)*t^2*p;                            #
y := (t^2-5*t-2+4/3*a)*t^2;                                        #
z := (2*t^2-(1+a)*t+4-8/3*a)*p;                                    #
F := factor(resultant(x*Z-z*X,y*Z-z*Y,t)/Z^4);                     #
factor(discrim(F,X));                                              #
factor(discrim(F,Y));                                              #
factor(discrim(F,Z));                                              #
\fi                                                                %

\smallskip\noindent
{\bf 36.} $A_7+2A_4+2A_2$.
$F=\Q$, $E=\Q(\sqrt{-3})$, $\omega^2+\omega+1=0$, $\{\omega,\omega/2\}\to2A_2$, $p=t$,
$$\split
x &= \omega t(2t^2 - 2\omega t + \omega^2)(4t^2 + 5t + 2),\\
y &= (2t^2 - t + 1)\big(2t^2 + (2+\omega)t - \omega^2\big)^2,\\
z &= \omega^2t^2(2t^2 + 2t + 1).
\endsplit$$
\if01                                                              %
alias(a=RootOf(Z^2+Z+1));                                          #
x := a*t*(2*t^2 - 2*a*t + a^2)*(4*t^2 + 5*t + 2);                  #
y := (2*t^2 - t + 1)*(2*t^2 + (2+a)*t - a^2)^2;                    #
z := a^2*t^2*(2*t^2 + 2*t + 1);                                    #
F := factor(resultant(x*Z-z*X,y*Z-z*Y,t)/Z^5);                     #
factor(discrim(F,X));                                              #
factor(discrim(F,Y));                                              #
factor(discrim(F,Z));                                              #
f := 32*Y^4                                                        #
 + (-16*X^2 - 288*X + 248)*Y^3                                     #
 + (2*X^4 + 96*X^3 + 570*X^2 - 2816*X - 198)*Y^2                   #
 + (-6*X^5 - 50*X^4 + 1130*X^3 +10650*X^2 + 1670*X - 238)*Y        #
 - 9*X^6 - 303*X^5 - 3430*X^4 - 13375*X^3 - 3530*X^2 + 973*X - 49; #
factor(subs(Z=1,F/f));                                             #
h := X^2+11*X-1;                                                   #
f1 := 32*Y^4                                                       #
 + (-16*X^2 - 288*X + 248)*Y^3                                     #
 + (2*X^4 + 96*X^3 +570*X^2 - 2816*X - 198)*Y^2                    #
 + (-6*X^3 + 16*X^2 + 948*X + 238)*h*Y - (9*X^2 + 105*X + 49)*h^2; #
factor(f1/f);                                                      #
\fi                                                                %

The implicit equation has rational coefficients:
$$
\split
  32y^4 &+ (-16x^2 - 288x + 248)y^3 + (2x^4 + 96x^3 +570x^2 - 2816x - 198)y^2\\
        & +(-6x^3 + 16x^2 + 948x + 238)hy - (9x^2 + 105x + 49)h^2;
\quad h = x^2+11x-1.
\endsplit
$$
the $A_4$ points are on the line $y=0$ at the roots of $h$.

\smallskip\noindent
{\bf 37.} $3A_6+A_1$, see [\refDegtMax]. $E=F=\Q$,
$[0,1,\infty]$,
$p=t^2-t+1$,
$$
x = (3t^2-3t+1)t^2,\qquad
y = (t^2+t+1)(t-1)^2,\qquad
z = (t^2-3t+3)t^2(t-1)^2;
$$
the curve is invariant under $(x:y:z)\mapsto(y:z:x)$, $t\mapsto 1/(1-t)$.
\if01                                                              %
x := (3*t^2-3*t+1)*t^2;                                            #
y := (t^2+t+1)*(t-1)^2;                                            #
z := (t^2-3*t+3)*t^2*(t-1)^2;                                      #
F := factor(resultant(x*Z-z*X,y*Z-z*Y,t)/Z^6);                     #
factor(discrim(F,X));                                              #
factor(discrim(F,Y));                                              #
factor(discrim(F,Z));                                              #
\fi                                                                %

\smallskip\noindent
{\bf 38.} $2A_6+A_4+A_2+A_1$.
$E=F=\Q(\sqrt{21})$, $a=\pm\sqrt{21}$,
[double roots of $z(t),\infty,0]$,
$p=2t^2+(3-a)t+6-{4\over3}a$,
$$\split
x &= \big( 2t^2 + (17+a)t - 2-4a \big)t^2,\\
y &= 10t^2 + (19-a)t + 10,\\
z &= \big( 2t^2 - (1+a)t + 8-2a \big)\big( 2t^2 + 6t + 1-a \big)^2.
\endsplit$$
\if01                                                              %
alias(a=RootOf(Z^2-21));                                           #
x := ( 2*t^2 + (17+a)*t - 2 - 4*a )*t^2;                           #
y :=  10*t^2 + (19-a)*t + 10 ;                                     #
z := (2*t^2 - (1+a)*t + 8-2*a)*(2*t^2 + 6*t + 1-a)^2;              #
F := factor(resultant(x*Z-z*X,y*Z-z*Y,t)/Z^6);                     #
factor(discrim(F,X));                                              #
factor(discrim(F,Y));                                              #
factor(discrim(F,Z));                                              #
\fi                                                                %

\smallskip\noindent
{\bf 39.} $A_6+A_5+2A_4$. $E=F=\Q(\sqrt{7})$,
%
%
$a=\pm\sqrt{7}$,
$[\infty$, roots of $q]$,
$p=t^2 - 21-8a$, $q=t^2+(6+a)t+7+3a$, $b=(7+2a)/3$,
$$
x = (t+2+a)(t+4+a)\,p,\quad
y = (t^2+12t+15+8a)\,p^2,\quad
z = \big(t^2-2at-b\big)\,q^2.
$$
\if01                                                              %
alias(a=RootOf(Z^2-7));  p := t^2 - 21 - 8*a;                      #
x := (t+2+a)*(t+4+a)*p;                                            #
y := (t^2+12*t+15+8*a)*p^2;                                        #
z := (t^2-2*a*t-(7+2*a)/3)*(t^2+(6+a)*t+7+3*a)^2;                  #
F := factor(resultant(x*Z-z*X,y*Z-z*Y,t)/Z^6);                     #
factor(discrim(F,X));                                              #
factor(discrim(F,Y));                                              #
factor(discrim(F,Z));                                              #
\fi                                                                %


\head \S2. The cases when $E\ne F$
\endhead

In this section we show that a curve $C$ of equisingularity types
{\bf 1}, {\bf 16}, {\bf 34}, {\bf 36} does not admit any parametrization
with coefficients in $F$. Recall that $F$ is the field of minimal degree
such that the curve is defined over $F$.
In cases {\bf 1} and {\bf 16} this is evident because $F$ is a real field but
the curve does not have any real smooth point.
In the other two cases we use a birational equivalence over $F$ of $C$ with
a plane conic $aX^2+bY^2=1$, $a,b\in F$. Such an equivalence exists due a classical result by Hilbert and Hurwitz [\refHH].
It follows that any rational curve defined by
an equation over $F$ admits a parametrization over a quadratic extension of $F$.
Moreover, a parametrization over $F$ exists if and only if the equation $aX^2+bY^2=1$
has a solution in $F$.

There are known several methods to find explicit formulas for a birational
equivalence of a rational curve with a plane conic defined over the same field.
The proof in [\refHH] provides a recursive
construction: at each step a curve of degree $d$ is mapped birationally to a curve
of degree $d-2$ using a generic $3$-dimensional system of adjoint
curves of degree $d-2$.
A more efficient algorithm based on similar ideas can be found in [\refSW].

Probably, the simplest modern proof of the Hilbert-Hurwitz theorem consists in saying
that the anticanonical linear system on the normalization of a rational curve
is defined over the same field; it has degree 2 and dimension 3, thus it embeds the curve to $\Bbb P^2$ as a conic. Of course, theoretically,
this proof can be developed up to an algorithm.
This was done in practice by van Hoeij in [\refvH].

\smallskip

Here we use a construction similar that in [\refACC, \refDegtACC]. We
consider a pencil of
cubic curves passing through 8 double points of $C$ including infinitely near points.
It happens that for curves considered here, the 8 points can be chosen so that
the pencil is defined over $F$ (note that this is not the case, for example,
for a generic rational sextic curve).
Each curve of this pencil has two intersections with $C$ outside the basepoints,
hence the pencil defines a double covering onto $\Bbb P^1$.
By choosing in addition a pencil of lines, we obtain a birational equivalence of
$C$ with a rational hyperelliptic curve which can be further transformed into a plane
conic. This method is illustrated below in the example of
the curves of equisingularity types {\bf 34} and {\bf 36}.
The underlying computations can be found in the files {\tt a7a4a4a2a2.mws} and
{\tt a7a6a6.mws} on the web page [\refOrevWebPg].

\medskip
{\bf 36.} $A_7+2A_4+2A_2$.
We start with the homogeneous defining equation $f(x,y,z)=0$
over $\Q$ given in \S1. The curve $C$ has an $A_7$ singularity at $(0:1:0)$
tangent to the axis $z=0$. It has two $A_4$ singularities on the axis $y=0$
at the roots of $h(x)=0$
and two $A_2$ singularities on the line $4x-y=0$.

Let $g_\lambda(x,y,z)=0$ be the pencil of cubic curves where $g_\lambda(x,y,1)$ is
$$
  8y^2 - \big(2x^2 +(4\lambda+24)x - 2\lambda+32\big)y + (\lambda x+7)h(x).
$$
It passes through all the four infinitely near points of $A_7$ and through each of the
other four singular points. 

Let us explain how we have found this expression for $g_\lambda$.
We consider a curve $g=0$ where
$g = 8y^2 - (2x^2 + c_1x + c_2)y + (c_3x+c_4)h(x)$ with indeterminate
coefficients $c_i$. Such a curve passes already through the both $A_4$
and through three infinitely near points at $A_7$.
Let $\gamma(t)=(t:1:\gamma_2 t^2+\gamma_3 t^3)$ where $\gamma_2=1/4$ and
$\gamma_3=c_1/32 - c_3/8$ are found from the condition $\ord_t(g(\gamma(t))=4$.
Then $\gamma$ defines a germ of a smooth curve at $(0:1:0)$ passing thought four
infinitely near points of the curve $g=0$.
Then the condition that $g$ passes through the fourth double point at $A_7$
reads $\ord_t f(\gamma(t))=8$. Due to the
chosen form of $f$ and $g$, this is equivalent to the vanishing of the
coefficient of $t^6$ in $f(\gamma(t))$ which is $(c_1-4c_3)^2/24$.

The two $A_2$ points are on the line $4x-y=0$, hence
the polynomial $f_L(t)=f(t,4t,1)$ has a factor $q(t)^2$, $\deg q=2$ (we have
$q(t)=\gcd(f_L,f'_L)=t^2-5t+7$; we can compute it
because $f$ is known).
So, we compute the remainder of the division of $g(t,4t,1)$ by $q(t)$
and we equate its coefficients to zero. By eliminating all variables but $c_3$
from these simultaneous equations, we obtain the above expression for $g_\lambda$.

From now on we pass to the affine chart $z\ne 0$. So, we set $z=1$ and we write
$f(x,y)$ and $g_\lambda(x,y)$ instead of $f(x,y,1)$ and $g_\lambda(x,y,1)$. Let
$$
    P(x,\lambda) = \Res_y( f(x,y), g_\lambda(x,y) ).
$$
It has the form $P=P_1(x,\lambda)P_2(x)$ where the factor $P_2(x)$ corresponds to the
base points of the pencil.
The curve $P_1(x,\lambda)=0$ is birationally equivalent to the curve $f(x,y)=0$.
The equivalence is given by the mapping which sends a generic point $(x,y)$
to $(x,\lambda)$ where $\lambda$ is found from the condition $g_\lambda(x,y)=0$ which is
a linear equation with respect to $\lambda$ (of course, the coordinate $y$
should be chosen generically enough to avoid a two-folded covering of
$f=0$ onto $P_1=0$). In our case we have
$$
\split
  P_1 &= (\lambda^4-15\lambda^3+63\lambda^2-486)x^2
                     + (4\lambda^4-78\lambda^3+234\lambda^2+2646\lambda-12636)x\\
  &\quad         - 2\lambda^4+156\lambda^3-3375\lambda^2+28053\lambda-80109,
 \qquad
  P_2 = 256\,q(x)^2h(x)^2.
\endsplit
$$
Let $D(\lambda) = \Discr_x P_1(\lambda,x)$.
Since the curve $P_1(x,\lambda)=0$ is rational, its projection onto the $\lambda$-axis has
exactly two branch points. They correspond to roots of $D$ of odd multiplicity.
The multiplicity of all other roots is even.
This means that $D=D_1D_2^2$ with $\deg_\lambda D_1=2$. It is easy to see that $P_1(x,\lambda)$
(and hence $f(x,y)$ as well) is  birationally equivalent 
the conic curve $D_1(\lambda)-u^2=0$ which in its turn is equivalent to
$Q(u,v)=0$ where $Q(u,v)=\Discr_\lambda(D_1(\lambda)-u^2)-v^2$.
In our case we have
$$
    D_1 = 6(\lambda^2-9\lambda+9), \quad D_2 = 2(\lambda-15)(\lambda^2-9\lambda+9),
    \quad Q = 24u^2+1620-v^2.
$$
Thus, we reduced the existence of a rational parametrization of
the initial curve to the existence of rational solutions to the equation
$6u^2+5w^2=1$ (note that $24=6\times2^2$ and $1620=5\times18^2$).
The Hilbert symbol $(6,5)_3$ is equal to $-1$, see [\refSerre; Th.~III-1],
hence there are no integer solutions.

\if01
g := Y^3 + (a*Z^2 + b*X*Z + 2*X^2)*Y + (28*Z^2 - 4*a*X*Z + (a-1)*X^2)*(c*Z + d*X);
e := expand(subs(X=x,Y=y,Z=z,g));

\fi

\medskip
{\bf 34.} $A_7+2A_6$. The proof is the same as
in case {\bf 36}. We find a pencil $g_\lambda$ passing through 4
double points at $A_7$ and through $2$ double points at each $A_6$.
Since we know a parametrization of the curve, the simplest way to do it
is to substitute the parametrization into
$y^2+(2x^2+c_1+c_2+x^2)y+(c_3x+c_4)h(x)$.
The result of the substitution is of the form $e(t)q(t)^2t^3$ where
$q(t)^2$ is the last factor of $y(t)$, see \S1-{\bf 34}.
Let $r(t)=e_0+e_1 t$ be the remainder of the division of $e(t)$ by $q(t)$.
Then $e(0)=e_0=e_1=0$ is the condition that $g_\lambda$ passes
through the required double points. So, we obtain
$$
\split
   g_\lambda(x,y) = y^2 &+ \big[2x^2 + 2(\lambda+a+9)x + (19a-7)/8\,\lambda + 49(a-1)/2\big]y\\
    &+ \big[\lambda x + (7-a)(\lambda+6)/2\big]\big(4x^2+4ax+7a-7\big).
\endsplit
$$
Further, we compute: $P(x,\lambda)=\Res_y(f,g_\lambda)/h^4$;
$\Discr_x P=D_1 D_2^2$ where
$$
\split
   &D_1 = (a-1)d(\lambda)/14,\qquad
   D_2 = 3^7(3a+1)(7\lambda+37a-35)d(\lambda)/1024\\
   &\text{and}\quad d(\lambda) = \lambda^2+(11a-1)\lambda-46a-54;\\
   &Q=\Discr_\lambda(D_1(\lambda)-u^2)-v^2=
(2u/a)^2\pi + \big(6(a-3)/7\big)^2 a - v^2
\endsplit
$$
where $\pi=(1-a)/2$ (recall that $a=\sqrt{-7}$).
Thus we obtain the equation
$$
    \pi X^2+aY^2 = 1.                                                 \eqno(\eqHilbS)
$$ 
It does not have a solution in $F$, since $(\pi,a)_\pi=-1$.
Since the ring of integers of $F$ is Euclidean, Hilbert symbols in $F$
are computed in the same way as in $\Q$. This means that the fact that 
(\eqHilbS) has no solution in $F$ admits the following ``high school algebra proof''.

Let $\Cal O_F$ be the ring of integers of $F$. If there is a solution to (\eqHilbS)
in $F$, then there exist $X,Y,Z\in\Cal O_F$ such that $\pi X^2+aY^2=Z^2$ and
there is no common non-unit divisor of $X$, $Y$, $Z$ in $\Cal O_F$.
Let us consider our equation reduced mod $\pi^3$.
The additive group $\Cal O_F$ is generated by 1 and $\pi$, hence the ideal $(\pi^3)$
is generated by $\pi^3$ and $\pi^4$. Since $\pi^4-3\pi^3=8$, we can chose
$\pi^3$ and $8$ as a base of $(\pi^3)$. Hence $\Cal O_F/(\pi^3)\cong\Z/8\Z$.
We have $\pi=-2-\pi^3\equiv 6\mod\pi^3$ and $a=1-2\pi\equiv1-2\times 6\equiv 5$.
Thus $(X,Y,Z)$ should be a solution to the congruence $6X^2+5Y^2\equiv Z^2\mod8$.
Since the only squares mod 8 are $0$, $1$, and $4$, we conclude that
$X$, $Y$, $Z$ represent even classes in $\Z/8\Z$. Since $\pi$ divides $2$
(indeed, $2=\pi-\pi^2$),
it follows that $X$, $Y$, $Z$ are all divisible by $\pi$ in $\Cal O_F$.


\head 3. How the parametrizations were found
\endhead

We used different approaches for different curves. Each approach is
discussed in a separate subsection of this section.
We used {\tt maple}, {\tt pari/gp}, {\tt Singular} and {\tt sage} software for computations.

\subhead 3.1. Parametric equations with indeterminate coefficients
\endsubhead
This approach works well when the most portion of $\mu(C)$ is contributed
by irreducible singularities. In all cases when we used this approach (except one case),
the computations were as easy as in [\refOrevFA] and they took
few seconds (sometimes minutes) of CPU time.

The most difficult case was
$A_7+2A_4+2A_2$ (no. {\bf 36}).
In this cases certain heuristic tricks were
applied to fasten the computations and despite them, many hours of parallel
work of many processors were needed to achieve the result.
An interesting coincidence is that Yang writes in [\refYang; p.~225] that:

\medskip
\hskip 0pt\hbox{\vbox{\pagewidth{120mm}\sl\noindent
``The largest discriminant is $3600$, the sextic curve is irreducible and its singularities
correspond to $A_7+2A_4+2A_2$. This is the only sextic curve whose discriminant reaches
$3600$.''}}

\medskip
I tried several ways to write a system of simultaneous equations for the
coefficients of the parametrization but every time the system, happened to
be so huge that I was not able to solve it on the available computers.
The system that I finally succeeded to solve is as follows.
We write a parametrization of a curve with indeterminate coefficients in the form
$$
   x=(1+a_1 t +a_2 t^2)t(t-1)^2,
\quad
   y=(1+b_1 t +b_2 t^2)t^2,
\quad
   z=(1+c_1 t +c_2 t^2)(t-1)^4.
$$
A generic curve $C$ of this form has an $A_3$ point at $(0:0:1)$
tangent to the line $y=0$ (corresponds to $t=0$ and $t=\infty$)
and an $A_4$ point at $(0:1:0)$ tangent to $z=0$ (corresponds to $t=1$).

Let $e(x,y,z)=z^2y - x^2z - \gamma x^3$ where
$\gamma=c_1 + b_1-2 a_1$ is found from the condition that $\ord_t e = 4$.
Then the curve $e=0$ is smooth at $(0:0:1)$ and it has tangency of order 4 with
the branch of $C$ at $t=0$.
We denote the coefficients of $t^{16}$ and $t^{15}$ in $e(x(t),y(t),z(t))$
by $e_1$ and $e_2$.
We have $e_1=c_2(b_2 c_2-a_2^2)$ and
$e_2=(2a_1-b_1-c_1)a_2^3 + (8c_2-c_1)a_2^2 + 2(b_2 c_1 - a_1 a_2)c_2 + (b_1-8b_2)c_2^2$.
Then the condition $e_1=e_2=0$ implies that the curve has an $A_k$ point
at $(0:1:0)$ with an odd $k\ge 7$ unless it is a multiple curve of a smaller degree.

\if01{ 
   x=(1+a1 t +a2 t^2)t(t-1)^2
   y=(1+b1 t +b2 t^2)t^2
   z=(1+c1 t +c2 t^2)(t-1)^4
   e=Expand[z^2 y - x^2 z - g x^3];
   Coefficient[e,t,0]
   Coefficient[e,t,1]
   Coefficient[e,t,2]
   e=e//.Solve[0==Coefficient[e,t,3],g][[1]];
   e1=Factor[Coefficient[e,t,15]]
   e2=Factor[Coefficient[e,t,16]]
}\fi 

Let $p(t)=t^3+at^2+bt+c$ where $a,b,c$ are indeterminates.
Let $e_3+e_4t+e_5t^2$ and $e_6+e_7 t+e_8 t^2$ be the remainders of the
division of
%
%
$(xy'-yx')/(t^2(t-1))$ and $(zy'-yz')/(t(t-1)^3)$ by $p(t)$.
Then the condition $e_3=\dots=e_8$ implies that
the curve has an $A_k$ point with an even $k\ge 2$ at each root of $p$.
Thus a generic solution to $e_1=\dots=e_8=0$ is a curve with $A_7+A_4+3A_2+A_1$.
It is hard but possible to eliminate all the variables except $b$ and $c$ from these
simultaneous equations. So we obtain an equation that we denote by $f[p]$
(in this notation we assume that $f$ is a polynomial in the coefficients of $p$,
not a polynomial in $p$).

Let $t\mapsto\varphi(t)$ be a parametrization of a $A_7+2A_4+2A_2$ curve such
that $0,\infty\mapsto A_7$ and $1\mapsto A_4$.
Let  $q(t) = t^2 + ut + v$ be the polynomial whose roots are mapped by $\varphi$ to
the $A_2$ points and let $\varphi(w)$ be the other $A_4$.
By changing the parameter of $\varphi$, we may obtain
three more parametrizations of the same curve such that $0,\infty\mapsto A_7$
and $1\mapsto A_4$, namely,
$$
\xalignat3
  &t\mapsto\varphi(wt)\;\;
     &&\text{maps:}\quad
        0,\infty\to A_7,\quad  1,w^{-1}\mapsto 2A_4,\quad 
        &&\text{roots of $q(wt)$}\mapsto 2A_2;
\\
  &t\mapsto\varphi(t^{-1})\;\;
     &&\text{maps:}\quad
        0,\infty\to A_7,\quad  1,w^{-1}\mapsto 2A_4,\quad
        &&\text{roots of $q(t^{-1})$}\mapsto 2A_2;
\\
  &t\mapsto\varphi(wt^{-1})\;\;
     &&\text{maps:}\quad
        0,\infty\to A_7,\quad  1,w\mapsto 2A_4,\quad
        &&\text{roots of $q(wt^{-1})$}\mapsto 2A_2.
\endxalignat
$$
Thus, $u,v,w$ satisfy the equations
$$
\split
  0 &=f\big[(t-w)(t^2+ut+v)\big]
  = f\big[(t-w^{-1})(t^2+uw^{-1}t+vw^{-2})\big]
\\
  &= f\big[(t-w^{-1})(t^2 + uw^{-1}t + w^{-1})\big]
  = f\big[(t-1)(t^2+uwv^{-1}t + w^2v^{-1})\big].
\endsplit
$$
By solving these simultaneous equations, we find the desired parametrization.


\subhead 3.2. Dual curve
\endsubhead
The degree of the dual curve $\check C$ is equal to $30-19-k$, $k=\#\Sing(C)$. For the set of singularities $(A_8+3A_2)+A_4+A_1$
(no. {\bf 33}) it is equal to $5$ (in all the other cases it is $\ge 6$) and the
set singularities of $\check C$ is $A_8 + A_4$. So,
we find a parametrization
$(\check x(t):\check y(t):\check z(t))$ of $\check C$ and we set
$$
   x = \check y' \check z - \check z' \check y, \quad
   y = \check z' \check x - \check x' \check z, \quad
   z = \check x' \check y - \check y' \check x.
$$

\smallskip\noindent
{\bf Remark.} One can check that the curves of equisingularity types with $k=5$,
i.~e., {\bf 26} $A_{10}+A_4+2A_2+A_1$, {\bf 36} $A_7+2A_4+2A_2$, {\bf 38}
 $2A_6+A_4+A_2+A_1$ are autodual.


\subhead 3.3. From an implicit equation to a parametrization
\endsubhead
In many cases we used the defining equations given in [\refACC] (for no. {\bf 1})
and in [\refDegtMax].
The problem of finding a parametrization of a rational curve defined by an implicit
equation is classical. In computer era it found new motivations in
geometric modeling and computer graphic. It worth to mention here that
Abhyankar and Bajaj [\refAB] gave an excellent presentation of a solution to
this problem (based on adjoint curves) accessible to people without any
algebro-geometric background. Their algorithm, however, is not appropiate
for our goals because the coefficients of the obtained parametrization do not
belong to a field extension of minimal degree.
 
We applied the method explained in \S2.
The only non-evident step is to find a solution in $F$
to the equation $Q=0$ (see \S2). This is an equation of the form $aX^2+bY^2=1$,
$a,b\in F$. In the most cases it was possible to choose the pencil $g_\lambda$ so
that the polynomial $D_1(\lambda)$ factorized into a product of two linear factors.
In the remaining cases, fortunately, the coefficients of $Q$ had the form $a=a_1a_2^2$, 
$b=b_1b_2^2$, $a_i,b_i\in F$, with $a_1$ and $b_1$ small enough that it was
possible to find a solution just by trying all pairs $(X,Y)$ one by one.

The most difficult case was $A_{10}+2A_4+A_1$ (no. {\bf 24}).
Starting with the implicit equation from [\refDegtMax],
performing the computations as in \S2 (with the pencil of cubics through $A_1$, $2A_4$,
and $5$ infinitely near points at $A_{10}$), and
changing the field generator to 
$\alpha$, $\alpha^3-\alpha^2-\alpha-1=0$ (found by {\tt pari}'s command
{\tt polred}),
we obtain very big expressions for $a$ and $b$. However, {\tt sage}'s commands
\smallskip
{\tt
\chardef\tthat`\^
K.<a>=NumberField(x\tthat 3-x\tthat 2-x-1)

K.ideal( {\it the coefficient} ).factor()
}
\smallskip\noindent
give a good hint how to extract quadratic factors of $a$ and $b$ and
we reduce our equation to
$(2-a)X^2-5a(a+2)Y^2-Z^2=0$.
Then we try all pairs $(X,Y)$ of the form
$X=x_0+x_1\alpha+x_2\alpha^2$,
$Y=y_0+y_1\alpha+y_2\alpha^2$, $x_i,y_i=0,\pm1,\pm2,\dots$
until the left hand side of the equation factorizes over $F$ into linear factors.
Thus we obtain a solution
$X=2+2a^2$, $Y=1+a-a^2$, $Z=2-a^2$.
The complete computation for $A_{10}+2A_4+A_1$ (from the equation in [\refDegtMax]
to the final formulas in \S2)
with detailed comments can be found in the file {\tt a10a4a4a1.mws} on the web page
[\refOrevWebPg].


\def\r{\ref}
\Refs

\r\no\refAB
\by S.S.~Abhyankar, C.L.~Bajaj
\paper Automatic parameterization of rational curves and surfaces III:
Algebraic plane curves
\jour Computer Aided Geometric Design \vol 5 \yr 1988 \pages 309--321
\endref

\r\no\refACC
\by E.~Artal Bartolo, J.~Carmona Ruber, J.I.~Cogolludo Agust{\'\i}n
\paper On sextic curves with big Milnor number
\inbook Trends in singularities \bookinfo Trends Math.
\publ Birkh\"auser \publaddr Basel \yr 2002 \pages 1--29
\endref

\r\no\refDegtPacif
\by A.~Degtyarev
\paper On plane sextics with double singular points
\jour Pacific J. Math. \vol 265 \yr 2013 \pages 327--348
\endref

\r\no\refDegtACC
\by A.~Degtyarev
\paper On the Artal--Carmona--Cogolludo construction
\jour J. Knot Theory Ramifications
\vol 23:5 \yr 2014
\endref

\r\no\refDegtMax
\by A.~Degtyarev
\paper Maximizing plane sextics with double points 
\jour Available online at\hbox to 20mm{}
http://www.fen.bilkent.edu.tr/$\tilde{\;}$degt/papers/equations.pdf
\endref

\r\no\refHH
\by D.~Hilbert, A.~Hurwitz
\paper \"Uber die Diophantischen Gleichungen vom Geschlecht Null
\jour Acta math. \vol 14 \yr 1890 \pages 217--224.
\endref

\r\no\refvH
\by M.~van Hoeij
\paper Rational parametrization of curves using canonical divisors
\jour J. Symbolic Computation \vol 23 \yr 1997 \pages 209--227
\endref

\r\no\refOkaPho
\by M.~Oka, D.T.~Pho
\paper Classification of sextics of torus type
\jour Tokyo J. Math. \vol 25 \yr 2002 \pages 399–433
\endref

\r\no\refOrevFA
\by S.Yu.~Orevkov
\paper A new affine $M$-sextic
\jour Funk. Anal. i Prilozhen.
\vol 32 \yr 1998 \issue 2 \pages 91--94
\transl English translation
\jour Funct. Anal. Appl. \vol 32 \yr 1998 \pages 141--143
\endref

\r\no\refOrevWebPg
\by S.Yu~Orevkov
\paper  Files related to this paper
\jour http://picard.ups-tlse.fr/$\tilde{\;}$orevkov/sextic-param.html
\endref

\r\no\refSW
\by J.R.~Sendra, F.~Winkler
\paper Parametrization of algebraic curves over optimal field extensions
\jour J. Symbolic Computation \vol 23 \yr 1997 \pages 191-207
\endref

\r\no\refSerre
\by J.-P.~Serre
\book Cours d'arithm\'etique
\publ Presses Universitaires de France \publaddr Paris \yr 1970
\endref

\r\no\refYang
\by J.G.~Yang
\paper  Sextic curves with simple singularities
\jour  Tohoku Math J. (2) \vol 48 \yr 1996 \pages 203--227
\endref

\endRefs

\enddocument